UDC 004.942

**Yakiv O. Kalinovsky,** Dr.Sc., Senior Researcher,
Institute for Information Recording National Academy of Science of Ukraine, Kyiv, Shpaka str. 2, 03113, Ukraine, E-mail: kalinovsky@i.ua

**Yuliya E. Boyarinova,** PhD, Associate Professor,
National Technical University of Ukraine "KPI", Kyiv, Peremogy av. 37, 03056, Ukraine,
E-mail: ub@ua.fm

**Alina S. Turenko**, postgraduate student,
Institute for Information Recording National Academy of Science of Ukraine, Kyiv, Shpaka str. 2, 03113, Ukraine, E-mail: asturenko@mail.ru

**Yana V. Khitsko,** Junior Researcher,
National Technical University of Ukraine "KPI", Kyiv, Peremogy av. 37, 03056, Ukraine,
E-mail: yannuary@yandex.ua

# Generalized quaternions and their relations with Grassmann-Clifford procedure of doubling

The class of non-commutative hypercomplex number systems (HNS) of 4-dimension, constructed by using of non-commutative Grassmann-Clifford procedure of doubling of 2-dimensional systems is investigated in the article and established here are their relationships with the generalized quaternions. Algorithms of performance of operations and methods of algebraic characteristics calculation in them, such as conjugation, normalization, a type of zero divisors are investigated. The considered arithmetic and algebraic operations and procedures in this class HNS allow to use these HNS in mathematical modeling.

**Key words:** quaternion, generalized quaternion, hypercomplex number system, zero divisor, pseudonorm, conjugation, Grassmann-Clifford procedure of doubling.

**Introduction**

There are many applications of hypercomplex number systems in science and technology. The system of quaternions, which is expansion of complex numbers, is of especially great importance. The system of quaternions, for the first time presented by V. R. Hamilton in work [1], found broad application in many scientific directions: in mechanics of a solid body, for the description of rotation in space, at the solution of problems of navigation, orientation and management of the movement; in computer animation, research of deformation of elastic designs, a filtration of color images, cryptography, etc.

**Definition of problem**

The main task of this work is in assignment of relations between the generalized quaternions and the systems of fourth dimension, constructed by means of non-commutative procedure of doubling of Grassmann-Clifford. Investigated are arithmetic and algebraic operations and procedures in this class of HNS that allows to use these HNS in mathematical modeling.

**The generalized quaternions**

Generalized quaternions were first employed in the description of a space-time isometry group by Kurt Gödel in his 1949 paper presenting his famous cosmological solution of the Einstein field equations [2]. These "Gödel quaternions" (also called "split quaternions" or "antiquaternions") belong to a real subalgebra of the complexified quaternion algebra which is not equivalent to the ordinary real quaternion algebra.

The generalized quaternions are also investigated in works of many other authors, for example [3 - 8]. We will analyze some results, which are received by them.

The generalized quaternion has a view:

$$q = a_1 e_1 + a_2 e_2 + a_3 e_3 + a_4 e_4 \tag{1}$$

where $a_i$ - are real numbers, and $e_i$ for $i = 2,...,4$ - imaginary units, which satisfy the following table of Cayley:

| $H_{\alpha\beta}$ | $e_1$ | $e_2$ | $e_3$ | $e_4$ |
|---|---|---|---|---|
| $e_1$ | $e_1$ | $e_2$ | $e_3$ | $e_4$ |
| $e_2$ | $e_2$ | $-\alpha e_1$ | $e_4$ | $-\alpha e_3$ |
| $e_3$ | $e_3$ | $-e_4$ | $-\beta e_1$ | $\beta e_2$ |
| $e_4$ | $e_4$ | $\alpha e_3$ | $-\beta e_2$ | $-\alpha\beta e_1$ |

$$\tag{2}$$

where $\alpha, \beta \in \mathbb{R}$.

As $e_1$ - is a real unit, than the generalized quaternion consists of two parts: valid and imaginary which are designated, respectively $S(q) = a_1 e_1$ and $V(q) = a_2 e_2 + a_3 e_3 + a_4 e_4$, thus it is possible to rewrite (1) in a format

$$q = S(q) + V(q) \tag{3}$$

When $\alpha = 1, \beta = 1$, the system of quaternions belongs to a class of the generalized quaternions. If to substitute these values of $\alpha$ and $\beta$ in (2), than we will receive Cayley table of quaternion system

| $H$ | $e_1$ | $e_2$ | $e_3$ | $e_4$ |
|---|---|---|---|---|
| $e_1$ | $e_1$ | $e_2$ | $e_3$ | $e_4$ |
| $e_2$ | $e_2$ | $-e_1$ | $e_4$ | $-e_3$ |
| $e_3$ | $e_3$ | $-e_4$ | $-e_1$ | $e_2$ |
| $e_4$ | $e_4$ | $e_3$ | $-e_2$ | $-e_1$ |

$$\tag{4}$$

Properties of the generalized quaternions are in detail considered in works [6, 7]. We will shortly give the main of them. Operation of multiplication is entered as well as for any other hypercomplex numbers. That is, carried out the rule:

$$\begin{aligned}(a_1 e_1 + a_2 e_2 + a_3 e_3 + a_4 e_4)(b_1 e_1 + b_2 e_2 + b_3 e_3 + b_4 e_4) = \\ = (a_1 b_1 - \alpha a_2 b_2 - \beta a_3 b_3 - \alpha\beta a_4 b_4) e_1 + \\ + (a_2 b_1 + a_1 b_2 - \beta a_4 b_3 + \beta a_3 b_4) e_2 + \\ + (a_3 b_1 + \alpha a_4 b_2 + a_1 b_3 - \alpha a_2 b_4) e_3 + \\ + (a_4 b_1 - a_3 b_2 + a_2 b_3 + a_1 b_4) e_4 \end{aligned} \tag{5}$$

Or by means of (3), it is possible to rewrite (5) in a view

$$qp = S(q)S(p) - \langle V(q), V(p) \rangle + S(q)V(p) + S(p)V(q) + V(p) \times V(q), \text{ where}$$

$$S(q) = a_1 e_1, \quad S(p) = b_1 e_1, \quad \langle V(q), V(p) \rangle = \alpha a_2 b_2 + \beta a_3 b_3 + \alpha\beta a_4 b_4,$$

$$V(p) \times V(q) = \beta(a_3 b_4 - a_4 b_3) e_2 + \alpha(a_4 b_2 - a_2 b_4) e_3 + (a_2 b_3 - a_3 b_2) e_4.$$

The conjugation for these numbers is entered as well as for quaternions, that is, if initial number is $q = a_1 e_1 + a_2 e_2 + a_3 e_3 + a_4 e_4$, than the conjugate number to it has a view:

$$\bar{q} = a_1 e_1 - a_2 e_2 - a_3 e_3 - a_4 e_4 \tag{6}$$

On the basis of conjugation, the norm is defined, which is determined from equality

$$N(q) = q\bar{q} = \bar{q}q = a_1^2 + \alpha a_2^2 + \beta a_3^2 + \alpha\beta a_4^2 \tag{7}$$

Also, in work [6] considered are concrete the HNS, depending on sign of $\alpha$ and $\beta$. The following cases are given:
1. $\alpha = \beta = 1$, then the hypercomplex number system $H_{\alpha\beta}$ is the system of quaternions;
2. $\alpha = 1$, $\beta = -1$, then the hypercomplex number system $H_{\alpha\beta}$ is the system of antiquaternions;
3. $\alpha = 1$, $\beta = 0$, then the hypercomplex number system $H_{\alpha\beta}$ is the system of semi-quaternions;
4. $\alpha = -1$, $\beta = 0$, then the hypercomplex number system $H_{\alpha\beta}$ is the system of semi-antiquaternions;
5. $\alpha = 0$, $\beta = 0$, then the hypercomplex number system $H_{\alpha\beta}$ is the system of $1/4$-quaternions.

**HNS of fourth dimension which are constructed by means of Grassmann-Clifford procedure of doubling**

As our researches showed [9], there are a relations between the systems received by means of non-commutative Grassmann-Clifford procedure of doubling and generalized quaternions. For establishment of this relationship, we will consider algorithm of construction and property of the systems, received by means of non-commutative Grassmann-Clifford procedure of doubling.

Investigated in work [9] the class of non-commutative HNS of fourth dimension consists of non-commutative doubling of HNS of second dimension by means of Grassmann-Clifford procedure of doubling. The basis of such HNS consists of four elements:

$$g = \{g_1, g_2, g_3, g_4\} = \{e_1 f_1, e_1 f_2, e_2 f_1, e_2 f_2\}.$$

Certainly, the two-character elements can be replaced by one-character elements with indexes. But it will be made later as at this stage it is appropriate to use the two-character elements.

Investigated class of HNS will be determined by the following conditions:
1. Elements of bases $e_1$ and $f_1$ - unit elements of their systems;
2. Elements of Cayley table of HNS g are products in form: $g_i g_k = e_j f_s e_t f_r$, whose values can be calculated by switching multipliers and using Cayley tables of HNS which are doubling; thus we will consider that $e_1$ and $f_1$ commutative with $e_2$ and $f_2$, that is $e_1 f_2 = f_2 e_1, e_2 f_1 = f_1 e_2$, and the last are non-commutative among themselves, that is $e_2 f_2 = -f_2 e_2$.

For example:
$$g_1 g_1 = e_1 f_1 e_1 f_1 = e_1 e_1 f_1 f_1 = e_1 f_1 = g_1;$$
$$g_2 g_1 = e_1 f_2 e_1 f_1 = e_1 e_1 f_2 f_1 = e_1 f_2 = g_2;$$
$$g_2 g_3 = e_1 f_2 e_2 f_1 = -e_1 e_2 f_2 f_1 = -e_2 f_2 = -g_4;$$
$$g_2 g_2 = e_1 f_2 e_1 f_2 = e_1 e_1 f_2 f_2 = e_1 f_2 f_2;$$
$$g_4 g_4 = e_2 f_2 e_2 f_2 = -e_2 e_2 f_2 f_2.$$

The last two examples of elements of the Cayley table can be finished only for concrete HNS which are doubling.

3. Elements of the Cayley table, which are under the main diagonal, but not in the first column, are opposite for elements, which are symmetric with respect to the main diagonal, that is

$$g_3 g_2 = -g_2 g_3;\ g_4 g_2 = -g_2 g_4;\ g_4 g_3 = -g_3 g_4.$$

Given these conditions, the generalized Cayley table for HNS of investigated class will be as follows:

|       | $g_1$ | $g_2$       | $g_3$        | $g_4$          |
|-------|-------|-------------|--------------|----------------|
| $g_1$ | $g_1$ | $g_2$       | $g_3$        | $g_4$          |
| $g_2$ | $g_2$ | $e_1 f_2 f_2$ | $-g_4$     | $-e_2 f_2 f_2$ |
| $g_3$ | $g_3$ | $g_4$       | $e_2 e_2 f_1$ | $e_2 e_2 f_2$ |
| $g_4$ | $g_4$ | $e_2 f_2 f_2$ | $-e_2 e_2 f_2$ | $-e_2 e_2 f_2 f_2$ |

(8)

As it is known [10-13] there are three isomorphisms classes of HNS of second dimension. We will choose from these classes one representatives: the system of complex numbers $C$, the system of double numbers $W$ and the system of dual numbers $D$.

As shown in the work [11], the first two operands in the operator of doubling can be commutative, as received Cayley tables differ only by order of lines and columns, that's why they are isomorphic.

Taking it into account, studied class of HNS consists of six representatives classes of isomorphism:
1. $D(C,C,4) = H$ – quaternion system;
2. $D(C,W,4) = AH$ – antiquaternion system;
3. $D(C,D,4) = D(C,D,4)$;
4. $D(W,W,4)$;
5. $D(D,D,4)$;
6. $D(W,D,4) = D(D,W,4)$;

Cayley tables of the above six classes of isomorphisms can be easily obtained by substituting in (8) the basic elements of the complex – $C$, double – $W$ and dual numbers – $D$, respectively

Having executed this algorithm for each of six cases, we will receive table 1

Table 1

**Cayley tables of hypercomplex number systems of fourth dimension**

| № | Designation | Cayley table |
|---|-------------|--------------|
| 1. | $H = D(C,C,4)$ | <table><tr><td>$H$</td><td>$e_1$</td><td>$e_2$</td><td>$e_3$</td><td>$e_4$</td></tr><tr><td>$e_1$</td><td>$e_1$</td><td>$e_2$</td><td>$e_3$</td><td>$e_4$</td></tr><tr><td>$e_2$</td><td>$e_2$</td><td>$-e_1$</td><td>$e_4$</td><td>$-e_3$</td></tr><tr><td>$e_3$</td><td>$e_3$</td><td>$-e_4$</td><td>$-e_1$</td><td>$e_2$</td></tr><tr><td>$e_4$</td><td>$e_4$</td><td>$e_3$</td><td>$-e_2$</td><td>$-e_1$</td></tr></table> |
| 2. | $AH = D(C,W,4)$ | <table><tr><td>$AH$</td><td>$e_1$</td><td>$e_2$</td><td>$e_3$</td><td>$e_4$</td></tr><tr><td>$e_1$</td><td>$e_1$</td><td>$e_2$</td><td>$e_3$</td><td>$e_4$</td></tr><tr><td>$e_2$</td><td>$e_2$</td><td>$-e_1$</td><td>$e_4$</td><td>$-e_3$</td></tr><tr><td>$e_3$</td><td>$e_3$</td><td>$-e_4$</td><td>$e_1$</td><td>$-e_2$</td></tr><tr><td>$e_4$</td><td>$e_4$</td><td>$e_3$</td><td>$e_2$</td><td>$e_1$</td></tr></table> |

| 3. | $D(C,D,4)$ | | $D(C,D,4)$ | $e_1$ | $e_2$ | $e_3$ | $e_4$ |
|---|---|---|---|---|---|---|---|
| | | | $e_1$ | $e_1$ | $e_2$ | $e_3$ | $e_4$ |
| | | | $e_2$ | $e_2$ | $-e_1$ | $e_4$ | $-e_3$ |
| | | | $e_3$ | $e_3$ | $-e_4$ | 0 | 0 |
| | | | $e_4$ | $e_4$ | $e_3$ | 0 | 0 |
| 4. | $D(W,W,4)$ | | $D(W,W,4)$ | $e_1$ | $e_2$ | $e_3$ | $e_4$ |
| | | | $e_1$ | $e_1$ | $e_2$ | $e_3$ | $e_4$ |
| | | | $e_2$ | $e_2$ | $e_1$ | $e_4$ | $e_3$ |
| | | | $e_3$ | $e_3$ | $-e_4$ | $-e_1$ | $e_2$ |
| | | | $e_4$ | $e_4$ | $-e_3$ | $-e_2$ | $e_1$ |
| 5. | $D(D,D,4)$ | | $D(D,D,4)$ | $e_1$ | $e_2$ | $e_3$ | $e_4$ |
| | | | $e_1$ | $e_1$ | $e_2$ | $e_3$ | $e_4$ |
| | | | $e_2$ | $e_2$ | 0 | $e_4$ | 0 |
| | | | $e_3$ | $e_3$ | $-e_4$ | 0 | 0 |
| | | | $e_4$ | $e_4$ | 0 | 0 | 0 |
| 6. | $D(W,D,4)$ | | $D(W,D,4)$ | $e_1$ | $e_2$ | $e_3$ | $e_4$ |
| | | | $e_1$ | $e_1$ | $e_2$ | $e_3$ | $e_4$ |
| | | | $e_2$ | $e_2$ | $e_1$ | $e_4$ | $e_3$ |
| | | | $e_3$ | $e_3$ | $-e_4$ | 0 | 0 |
| | | | $e_4$ | $e_4$ | $-e_3$ | 0 | 0 |

**Relations between the systems received by means of Grassmann-Clifford procedure of doubling and the generalized quaternions**

As we can see from table 1, multiplication of basic elements of system of quaternions satisfies (4).

If to analyze the received results further, it is possible to see that the system of antiquaternions $AH$ corresponds to the second case of the generalized quaternions for $\alpha=1$, $\beta=-1$, and is called in work [6] as a split-quaternion system.

The system $D(C,D,4)$ corresponds to the third case of the generalized quaternions – the system of semi-quaternion for $\alpha=1$, $\beta=0$.

The system $D(W,D,4)$ corresponds to the fourth case - the system of semi-split quaternion for $\alpha=-1$, $\beta=0$.

And, at last, fifth case $1/4$ - quaternion for $\alpha=0$, $\beta=0$ corresponds to the system $D(D,D,4)$ from the table 1.

In case of the generalized quaternions are considered fifth separate cases, and by means of Grassmann-Clifford procedure of doubling - are received six systems. That is, the system $D(W,W,4)$ for the generalized quaternions isn't considered. Having analyzed Cayley table of this system, authors of this article defined that it corresponds to the multiplication table of basic elements of the generalized quaternions (2) for a case of $\alpha=-1$, $\beta=-1$.

**Research of properties of the received class of HNS**

We will designate numbers of each of these systems as $w=a_1e_1+a_2e_2+a_3e_3+a_4e_4$, where: $a_i\in R$.

In these systems operation of multiplication is carried out by the rule of multiplication of any two hypercomplex numbers, in respect that Cayley table of the considered system. Or, substituting in (5) for each of systems the values $\alpha$ and $\beta$. In both cases we will receive identical results.

In the work [14], the norm of hypercomplex number generally is determined by a formula

$$N(w) = \sum_{i=1}^{n} \gamma_{ij}^{k} a_i,$$

where $\gamma_{ij}^{k}$ - structural constants of hypercomplex number system, which are defined from table 1. On this basis the norm of matrix is constructed [14], having calculated the determinant of which, we will receive the norm of hypercomplex number.

By analogy to the theory of quaternions, we will call a root of the norm a pseudonorm of hypercomplex number, which will be denoted as $N(w)$:

It should be noted that for each of these systems the matrix of norm will have other appearance, and according to it, it will have representations of a pseudonorm, that it is possible to see from table 2.

Table 2

**Pseudonorm**

| № | HNS | Pseudonorm |
|---|---|---|
| 1. | $H$ | $N(w) = a_1^2 + a_2^2 + a_3^2 + a_4^2$ |
| 2. | $AH$ | $N(w) = a_1^2 + a_2^2 - a_3^2 - a_4^2$ |
| 3. | $D(C,D,4)$ | $N(w) = a_1^2 + a_2^2$ |
| 4. | $D(W,W,4)$ | $N(w) = a_1^2 - a_2^2 - a_3^2 + a_4^2$ |
| 5. | $D(D,D,4)$ | $N(w) = a_1^2$ |
| 6. | $D(W,D,4)$ | $N(w) = a_1^2 - a_2^2$ |

If in equality (7) for each of systems, to substitute the values $\alpha$ and $\beta$, we will receive the same results. Apparently from table 2, in some systems the pseudonorm can be negative. It is possible to show that the pseudonorm entered by such method is multiplicative for each of considered systems, i.e. equality is carried out:

$$N(w_1 w_2) = N(w_1) N(w_2). \qquad (9)$$

As is offered in [14], definition of conjugate is entered on the basis of equality

$$w\overline{w} = N(w),$$

where $\overline{w} = b_1 e_1 + b_2 e_2 + b_3 e_3 + b_4 e_4$.

Substituting introduced notation and using Table 2, and equating coefficients of the same basic elements we obtain a linear algebraic system with respect to variables $b_1, b_2, b_3, b_4$.

For hypercomplex number system $D(D,D,4)$, for example, this linear algebraic system is

$$\begin{cases} a_1 b_1 = a_1^2 \\ a_1 b_2 + a_2 b_1 = 0 \\ a_1 b_3 + a_3 b_1 = 0 \\ a_1 b_4 + a_4 b_1 + a_2 b_3 - a_3 b_2 = 0 \end{cases}, \qquad (10)$$

which solutions have the form: $b_1 = a_1, b_2 = -a_2, b_3 = -a_3, b_4 = -a_4$.

Therefore, if the original number $w = a_1 e_1 + a_2 e_2 + a_3 e_3 + a_4 e_4$, that the conjugate number to it has view:

$$\overline{w} = a_1 e_1 - a_2 e_2 - a_3 e_3 - a_4 e_4. \tag{11}$$

In spite of the fact that for each of the systems given in table 1, the linear algebraic system (10) will have other appearance, representation of the interfaced number for each case will have an appearance (11), that coincides with (6).

Also, for the considered class, the signs of zero divisors are defined.

Not equal to zero hypercomplex number $w_1 \neq 0$ is called as zero divider if there is such hypercomplex number $w_2 \neq 0$, that their product is equal to zero $w_1 w_2 = 0$, and it means the same ratio between their pseudonorms $N(w_1 w_2) = 0$.

On the basis of (9), the pseudonorm of zero divisor is equal to zero

$$N(w_1) = 0. \tag{12}$$

From (12) follows the signs of zero divisor in any considered hypercomplex number system, (except system of quaternions for which, according to Frobenius's theorem, there are no zero divisors), which we will give in table 3

Table 3

**Signs of zero divisor**

| № | HNS | Sign of zero divisor |
|---|---|---|
| 1. | $H$ | - |
| 2. | $AH$ | $a_1^2 + a_2^2 = a_3^2 + a_4^2$ |
| 3. | $D(C, D, 4)$ | $a_1^2 + a_2^2 = 0$ |
| 4. | $D(W, W, 4)$ | $a_1^2 + a_4^2 = a_3^2 + a_2^2$ |
| 5. | $D(D, D, 4)$ | $a_1^2 = 0$ |
| 6. | $D(W, D, 4)$ | $a_1^2 = a_2^2$ |

**Conclusions**

In this work, established are the relations between separate cases of the generalized quaternions and the systems of fourth dimension, constructed by means of non-commutative Grassmann-Clifford procedure of doubling of the systems of second dimension. Investigated are their arithmetic and algebraic properties that allows to draw a conclusion on possibility of their usage for creation of various mathematical models.